\documentclass{article}
\usepackage{amsmath,amssymb}
\usepackage{graphicx} 
\usepackage{blkarray}
\usepackage{xcolor}
\usepackage{hyperref}
\usepackage{subfig}
\usepackage{comment}

\usepackage{array}

\usepackage{lineno}

\title{Relations between the properties of a complete rooted tree and the properties of a distribution of lengths of randomly generated strings}
\author{Yurii Lahodiuk \\ \texttt{yura.lagodiuk@gmail.com}}
\date{}

\begin{document}

\maketitle

\begin{abstract}
In scope of this paper we show that there are relations between the properties of a complete $m$-ary rooted tree of height $n$ and the expectation and variance of the distribution of lengths of random strings generated from the alphabet $\{ \alpha_1, \alpha_2, \dots \alpha_m \}$ until we see $n$ instances of a specific symbol in a row. Based on this we demonstrate a new interpretation for the integer sequence A286778.
\end{abstract}


\subsection*{Overview}

Consider a complete $m$-ary rooted tree graph of height $n$: $G_{m,n}=(V, E)$.
Let's denote the total number of edges in such tree as $T_{m,n} := |E| = m \cdot {m^n - 1 \over m - 1}$.
For every node $v \in V$ there always exists a unique path $\pi(v)$ from the node $v$ to the root node (we treat paths as sets of edges).
For any pair of nodes $(a,b) \in V \times V$ we can find the length of the common sub-path to the root node: $|\pi(a) \cap \pi(b)|$.
Let's denote the sum of the common sub-path lengths over all 2-tuples of nodes of such tree as $S_{m,n} := \sum_{(a,b) \in V \times V} |\pi(a) \cap \pi(b)|$.
\begin{figure}[!hbt]
    \centering
    \subfloat[$T_{2,2}=6$ and $S_{2,2}=22$.] 
    {{ \includegraphics[scale=0.5]{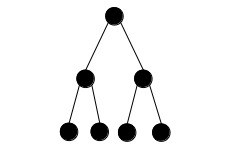} \label{fig:intro_m_2_figure} }}
    \qquad
    \subfloat[$T_{3,2}=12$ and $S_{3,2} = 57$] 
    {{ \includegraphics[scale=0.5]{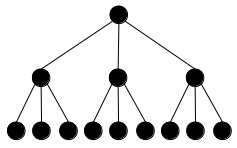} \label{fig:intro_m_3_figure} }}
    \caption{Examples of the trees and the corresponding values of $T_{m,n}$ and $S_{m,n}$.}
\end{figure}

On the other hand, consider a process that generates a random string as follows: starting from the empty string we pick a random symbol from the alphabet $\{ \alpha_1, \alpha_2, \dots \alpha_m \}$ and append it to the string, the process continues until we see $n$ instances of a specific symbol (say, $\alpha_1$) in a row. Let's introduce a random variable $\xi_{m,n}$ that represents a length of the string generated by the described process.

We will prove that the following relations between $T_{m,n}$ and the expectation of $\xi_{m,n}$, and between $S_{m,n}$ and the variance of $\xi_{m,n}$, are true for all $m,n \geq 1$:
\begin{equation}
\label{variance_proportional_to_sum_of_paths_length}
    \mathrm{Var}[\xi_{m,n}] = (m-1) \cdot S_{m,n}
\end{equation}
\begin{equation}
\label{expectation_proportional_to_tree_size}
        \mathbb{E}[\xi_{m,n}] = T_{m,n}
\end{equation}

In the table \ref{table_with_examples} you can see few examples:
\begin{table}[!h]
\centering
\begin{tabular}{|c|c|c|c|c|}
\hline
(m,n) & $T_{m,n}$ & $\mathbb{E}[\xi_{m,n}]$ & $S_{m,n}$ & $\mathrm{Var}[\xi_{m,n}]$ \\ \hline
(2,2) &    6       &      6        &     22      &   22  \\ \hline
(2,3) &    14       &      14        &     142      &   142  \\ \hline
(2,4) &    30       &      30        &     734      &   734  \\ \hline
(2,5) &    62       &      62        &     3390      &   3390  \\ \hline
\hline
(3,2) &     12      &      12        &     57      &   114 \\ \hline
(3,3) &      39     &       39       &      678      &   1356 \\ \hline
(3,4) &      120     &       120       &      6834      &   13668 \\ \hline
\hline
(4,2) &     20      &      20        &     116      &   348 \\ \hline
(4,3) &      84     &       84       &      2228      &   6684 \\ \hline
(4,4) &      340     &       340       &      37812      &   113436 \\ \hline
\end{tabular}
\caption{Note, that the relations between the items in the rows obey the equations \eqref{variance_proportional_to_sum_of_paths_length} and \eqref{expectation_proportional_to_tree_size}.}
\label{table_with_examples}
\end{table}

For proving the relations \eqref{variance_proportional_to_sum_of_paths_length} and \eqref{expectation_proportional_to_tree_size} in case when $m \geq 2$ we will use the closed-form expressions for $T_{m,n}$, $\mathbb{E}[\xi_{m,n}]$, $S_{m,n}$ and $\mathrm{Var}[\xi_{m,n}]$. 

The computation of expectation $\mathbb{E}[\xi_{m,n}]$ and variance $\mathrm{Var}[\xi_{m,n}]$ is described in the literature, for instance: the approach based on generating functions is described in \cite{graham_knuth_patashnik_concrete_mathematics}, and the approach based on martingales is described in \cite{probability_with_martingales}, also the closed-form expression for $\mathbb{E}[\xi_{m,n}]$ is provided in \cite{nielsen_expectation}. 

For completeness, in scope of these notes we will derive the closed-form expressions for $\mathbb{E}[\xi_{m,n}]$ and $\mathrm{Var}[\xi_{m,n}]$. We will use the transfer matrix method \cite{transfer_matrix_method}, which is similar to the approach based on the absorbing Markov chain. In the section \ref{expectation_and_variance} we will construct a matrix $W$ based on the adjacency matrix of a digraph that corresponds to the described string-generation problem. Using the fact, that the spectral radius of $W$ is less than $1$ (which we will prove in the section \ref{infinite_series_convergence}) we will derive the matrix-form expressions for computing the expectation and variance (where $\mathbb{I}$ is an identity matrix):
\begin{equation} 
\label{expectation_matrix_equation}
\begin{split}
	\mathbb{E}[\xi_{m,n}] = (0, \dots , 0, 1) \cdot W \cdot (\mathbb{I} - W)^{-2} \cdot (1, 0, \dots , 0)^\mathsf{T}
\end{split}    
\end{equation}
\begin{equation} 
\label{variance_matrix_equation}
\begin{split}
    \mathrm{Var}[\xi_{m,n}]  = (0, \dots , 0, 1) \cdot W \cdot (\mathbb{I} + W) \cdot (\mathbb{I} - W)^{-3} \cdot (1, 0, \dots , 0)^\mathsf{T} - \mathbb{E}[\xi_{m,n}]^2
\end{split}    
\end{equation}

Afterwards (in the section \ref{closed_form_expectation_and_variance}), taking into account relations between the entries of the matrices $W$, $(\mathbb{I} + W)$, and $(\mathbb{I} - W)^{-1}$, we will derive from the equations \eqref{expectation_matrix_equation} and \eqref{variance_matrix_equation} the closed-form expressions for expectation and variance:
\begin{equation} 
\label{coin_tossing_expectation_initial_equation}
\begin{split}
	\mathbb{E}[\xi_{m,n}] = m \cdot {m^n - 1 \over m - 1}
\end{split}    
\end{equation}
\begin{equation}
\label{coin_tossing_variance_initial_equation}
\begin{split}
    \mathrm{Var}[\xi_{m,n}] & = {m \over (m-1)^2} \cdot (m^{2n+1} - (2n+1) \cdot m^{n+1} + (2n+1) \cdot m^n - 1)
\end{split}
\end{equation}

In the section \ref{paths_in_tree_section} we will derive the closed-form expression for $S_{m,n}$ using the combinatorial counting techniques:
\begin{equation}
\label{sum_path_length_initial_equation}
\begin{split}
    S_{m,n} & = {m \over (m-1)^3} \cdot (m^{2n+1} - (2n+1) \cdot m^{n+1} + (2n+1) \cdot m^n - 1)
\end{split}
\end{equation}

The equation \eqref{variance_proportional_to_sum_of_paths_length} follows from the equations \eqref{coin_tossing_variance_initial_equation} and \eqref{sum_path_length_initial_equation}. 
And the equation \eqref{expectation_proportional_to_tree_size} follows from the equation \eqref{coin_tossing_expectation_initial_equation} and the definition of $T_{m,n}$.

In case when $m=2$ we have $\mathbb{E}[\xi_{2,n}] = T_{2,n}$ and $\mathrm{Var}[\xi_{2,n}] = S_{2,n}$ for all $n \geq 1$. 
While it is known that both $\mathbb{E}[\xi_{2,n}]$ and $T_{2,n}$ are described by the integer sequence A000918 \cite{oeis_A000918}, and it is known that $S_{2,n}$ is described by the integer sequence A286778 \cite{oeis_A286778}, we obtain a new interpretation for the integer sequence A286778: \textit{this sequence describes $\mathrm{Var}[\xi_{2,n}]$ -- a variance of the number of tosses of a fair coin until we see $n$ heads in a row}. In this case we can derive from the equation \eqref{coin_tossing_variance_initial_equation} the same formula as provided in \cite{oeis_A286778} for computing the values of A286778:
\begin{equation} 
\label{oeis_A286778_formula}
    a_n=4 \cdot 2^{2n} - (4n+2) \cdot 2^n - 2
\end{equation}

As a closing remark, it worth to admit, that it might be interesting to investigate further: whether it is possible to find similar relations for the distributions of lengths of randomly generated strings, that are generated until we observe any particular (given in advance) word for the given alphabet.


\section{Expectation and variance of the lengths of generated strings}
\label{expectation_and_variance}

In case if $m=1$ (the alphabet contains only one symbol) the expressions for expectation and variance are trivial: $\mathbb{E}[\xi_{1,n}]=n$ and $\mathrm{Var}[\xi_{1,n}]=0$ (and the equations \eqref{variance_proportional_to_sum_of_paths_length} and \eqref{expectation_proportional_to_tree_size} are true in this case), so in the remaining part of the paper we will focus on the case when $m \geq 2$.

\begin{figure}[!hbt]
    \includegraphics[scale=0.4]{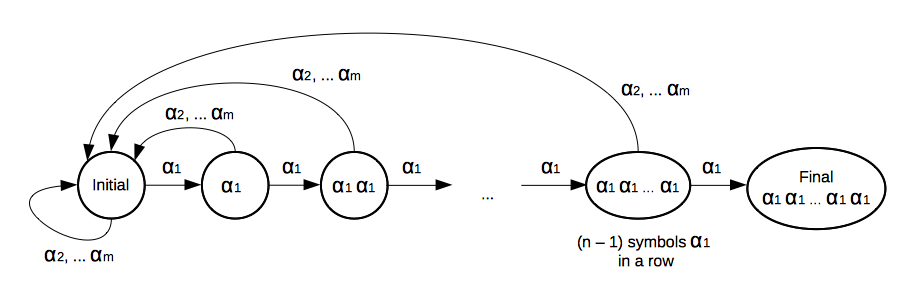}
    \caption{After choosing the random symbol from the alphabet we follow the edge, which contains the chosen symbol (we start from the ``initial" vertex).}
    \label{fig:walk_in_graph}
\end{figure}

We can represent our string generation problem as a walk in the digraph displayed in the figure \ref{fig:walk_in_graph} (each edge is associated with a subset of symbols from our alphabet: either $\{ \alpha_2, \dots \alpha_m \}$, or $\{ \alpha_1 \}$).   

This graph has a $(n+1) \times (n+1)$ adjacency matrix \eqref{eq:matrix_A_structure}, where the outgoing edges of vertices are represented as column-vectors, and the values of the items of this matrix represent the sizes of the subsets of symbols associated with the corresponding edges.

\begin{equation}
\label{eq:matrix_A_structure}
A=
\begin{blockarray}{ccccccc}
Initial & \alpha_1 & \alpha_1 \alpha_1 & \dots & \alpha_1 \alpha_1...\alpha_1 & Final  \\
\begin{block}{(cccccc)c}
\\
(m-1)   & (m-1) & (m-1)  & \dots    & (m-1)  & 0         & Initial      \\ \\
1         & 0        & \dots   & \dots    & 0         & 0         & \alpha_1      \\ \\
0         & 1        & 0         & \dots    & \vdots & \vdots & \alpha_1 \alpha_1 \\ \\
\vdots & \ddots & \ddots & \ddots  & \vdots & \vdots & \vdots \\ \\
\vdots &            & \ddots & \ddots  & 0         & \vdots & \alpha_1 \alpha_1 ... \alpha_1 \\ \\
0        & \dots    & \dots   & 0          & 1         & 0        & Final      \\
\\
\end{block}
\end{blockarray}
\end{equation}

We are interested in the number of walks, that arrive into the final vertex after $k$ steps (starting from the initial vertex). The matrix $A^k$ allows to obtain the number of walks of length $k$ between any pair of vertices. So, the number of walks between the initial vertex and between the final vertex is a leftmost bottom element of $A^k$:
\begin{equation}
    (0, \dots , 0, 1) \cdot A^k \cdot (1, 0, \dots , 0)^\mathsf{T}
\end{equation}

Given the alphabet size $m$ the total number of walks of length $k$ is $m^k$, therefore the probability that the random walk of the length $k$ will finish in the final vertex is:
\begin{equation}
\begin{split}
    p_k & = {1 \over m^k} \cdot (0, \dots , 0, 1) \cdot A^k \cdot (1, 0, \dots , 0)^\mathsf{T} \\
           & =  (0, \dots , 0, 1) \cdot \left( {1 \over m} \cdot A \right)^k \cdot (1, 0, \dots , 0)^\mathsf{T}
\end{split}
\end{equation}

Let's denote $W:={1 \over m} \cdot A$. Also, let's introduce a discrete random variable $\xi_{m,n}$, that represents a length of the string, generated until we observed $n$ symbols $\alpha_1$ in a row ($\xi_{m,n}$ takes values from the set $\{0, 1, 2, 3, \dots \}$). Thus, the expected length of the string is:
\begin{equation} 
\label{eq:expected_amount_of_tosses_infinite_sum}
    \mathbb{E}[\xi_{m,n}] = \sum_{k=0}^{\infty} k \cdot p_k = (0, \dots , 0, 1) \cdot \left( \sum_{k=0}^{\infty} k \cdot W^k \right) \cdot (1, 0, \dots , 0)^\mathsf{T}
\end{equation}

In the section \ref{infinite_series_convergence} we will investigate the spectral radius of $W$, and we will show that $\rho(W) < 1$. Thus, we deal with the Neumann series: $\sum_{k=0}^{\infty} W^k=(\mathbb{I} - W)^{-1}$. Consequently, we can transform the infinite sum in equation \eqref{eq:expected_amount_of_tosses_infinite_sum} using the same technique as with formal power series: 
\begin{equation}
    \sum_{k=0}^{\infty} k \cdot x^k = x \cdot {d \over dx} \left( \sum_{k=0}^{\infty} x^k \right) = {x \over (1-x)^2}
\end{equation}

Hence, the expectation can be expressed as:
\begin{equation}
\label{eq:expected_amount_of_tosses_via_inverted_matrix}
    \mathbb{E}[\xi_{m,n}] = (0, \dots , 0, 1) \cdot W \cdot \left( (\mathbb{I} - W)^{-1} \right)^2 \cdot (1, 0, \dots , 0)^\mathsf{T}
\end{equation}

The product of matrices $W \cdot (\mathbb{I} - W)^{-1}$ is commutative (the proof is in the section \ref{matrix_products_properties_section}), so the order of matrices in the equation \eqref{eq:expected_amount_of_tosses_via_inverted_matrix} does not matter.

Now, let's compute the variance:
\begin{equation}
    \mathrm{Var}[\xi_{m,n}] = \mathbb{E}[\xi_{m,n}^2] - \mathbb{E}[\xi_{m,n}]^2
\end{equation}

Let's compute $\mathbb{E}[\xi_{m,n}^2]$:
\begin{equation}
    \mathbb{E}[\xi_{m,n}^2] = \sum_{k=0}^{\infty} k^2 \cdot p_k = (0, \dots , 0, 1) \cdot \left( \sum_{k=0}^{\infty} k^2 \cdot W^k \right) \cdot (1, 0, \dots , 0)^\mathsf{T}
\end{equation}

As far as $\rho(W) < 1$, we can transform the sum $\sum_{k=0}^{\infty} k^2 \cdot W^k$ using the same technique as with formal power series: 
\begin{equation}
    \sum_{k=0}^{\infty} k^2 \cdot x^k = x \cdot {d \over dx} \left( \sum_{k=0}^{\infty} k \cdot x^k \right) = {x \cdot (1 + x) \over (1-x)^3}
\end{equation}

Hence, $\mathbb{E}[\xi_{m,n}^2]$ could be expressed as:
\begin{equation}
\label{eq:expectation_n_2}
    \mathbb{E}[\xi_{m,n}^2] = (0, \dots , 0, 1) \cdot W \cdot (\mathbb{I} + W) \cdot \left( (\mathbb{I} - W)^{-1} \right)^3 \cdot (1, 0, \dots , 0)^\mathsf{T}
\end{equation}

The products of matrices $W \cdot (\mathbb{I} + W)$ and $(\mathbb{I} + W) \cdot (\mathbb{I} - W)^{-1}$ are commutative (the proof is in the section \ref{matrix_products_properties_section}), so the order of matrices in the equation \eqref{eq:expectation_n_2} does not matter.

From equations \eqref{eq:expected_amount_of_tosses_via_inverted_matrix} and \eqref{eq:expectation_n_2} we obtain the matrix-form expression for the variance:
\begin{equation}
\label{eq:matrix_form_expression_for_variance}
\begin{split}
    \mathrm{Var}[\xi_{m,n}]  = & \mathbb{E}[\xi_{m,n}^2] - \mathbb{E}[\xi_{m,n}]^2 \\
    				   = & (0, \dots , 0, 1) \cdot W \cdot (\mathbb{I} + W) \cdot \left( (\mathbb{I} - W)^{-1} \right)^3 \cdot (1, 0, \dots , 0)^\mathsf{T} - \\
				      & - \left( (0, \dots , 0, 1) \cdot W \cdot \left( (\mathbb{I} - W)^{-1} \right)^2 \cdot (1, 0, \dots , 0)^\mathsf{T} \right)^2 \\
\end{split}
\end{equation}

Generally speaking, for the cases when we generate random strings until we encounter any arbitrary word constructed from the given alphabet (not only until we see $n$ instances of a specific symbol in a row) -- the matrix-form expressions for expectation and variance \eqref{eq:expected_amount_of_tosses_via_inverted_matrix} and \eqref{eq:matrix_form_expression_for_variance} are true as long as the series, that involve the corresponding matrix $W$, converges. But, still, in scope of this paper we focus on the random string generation process that continues until we see $n$ instances of a specific symbol in a row.

\section{The closed-form expressions for expectation and variance}
\label{closed_form_expectation_and_variance}

Taking into account the structure of the matrix $A$ from the equation \eqref{eq:matrix_A_structure} and that $W:={1 \over m} \cdot A$, we can notice that the values of entries of the matrix $(\mathbb{I} - W)^{-1}$ obey a simple pattern:
\begin{equation}
\label{eq:one_minus_w_inv}
\begin{blockarray}{cccccccc}
\begin{block}{(cccccccc)}
\\
m^n       &(m^n - m^1) &(m^{n} - m^2)             & \dots    &(m^n - m^k)              & \dots   & (m^n - m^n) \\
\vdots    & m^{n-1}       & (m^{n-1} - m^1)        & \dots     &(m^{n-1} - m^{k-1})  & \dots   &  \vdots          \\ 
\vdots    & \vdots          & m^{n-2}                    & \dots     &\vdots                      & \dots   &  \vdots          \\ 
\vdots    & \vdots          & \vdots                       &             & \vdots                      &            &  \vdots          \\ 
\vdots    & \vdots          & \vdots                       & \ddots  & (m^{n-k+1}-m^1)     &            &  \vdots          \\ 
\vdots    & \vdots          & \vdots                       &             & m^{n-k}                   &             &  \vdots          \\ 
\vdots    & \vdots          & \vdots                       &             & \vdots                      &            &  \vdots          \\ 
m^2       & m^2             & m^2                          & \dots    & m^2                        &             & (m^2 - m^2)  \\
m^1       & m^1             & m^1                          & \dots    & m^1                        & \ddots  & (m^1 - m^1)  \\
m^0       & m^0             & m^0                          & \dots    & m^0                        & \dots    & m^0              \\
\\
\end{block}
\end{blockarray}
\end{equation}

Considering the structure of the matrix $W$ we can compute $(0, \dots, 0, 1) \cdot W$:
\begin{equation}
\label{equation_0_0_1_times_w}
    (0, \dots, 0,  0, 1) \cdot W = \left( 0, \dots, 0, {1 \over m}, 0 \right)
\end{equation}

Let's compute $(\mathbb{I} - W)^{-1} \cdot (1, 0, \dots , 0)^\mathsf{T}$ using the matrix from \eqref{eq:one_minus_w_inv}:
\begin{equation}
\label{equation_i_sub_w_inv_times_0_0_1_T}
    (\mathbb{I} - W)^{-1} \cdot (1, 0, \dots , 0)^\mathsf{T} = (m^n, m^{n-1}, m^{n-2}, \dots, m^0)^\mathsf{T}
\end{equation}

Using the similar vector-matrix multiplications as in \eqref{equation_0_0_1_times_w} and \eqref{equation_i_sub_w_inv_times_0_0_1_T} we can compute the closed-form expression for $\mathbb{E}[\xi_{m,n}]$:
\begin{equation}
\label{equation_e_xi}
\begin{split}
        \mathbb{E}[\xi_{m,n}] & = (0, \dots , 0, 1) \cdot W \cdot \left( (\mathbb{I} - W)^{-1} \right)^2 \cdot (1, 0, \dots , 0)^\mathsf{T} \\
                                           & = \left( 0, \dots, 0, {1 \over m}, 0 \right) \cdot (\mathbb{I} - W)^{-1} \cdot (m^n, m^{n-1}, m^{n-2}, \dots, m^0)^\mathsf{T} \\
                                           & = \left( \left( 0, \dots, 0, {1 \over m}, 0 \right) \cdot (\mathbb{I} - W)^{-1} \right) \cdot (m^n, m^{n-1}, m^{n-2}, \dots, m^0)^\mathsf{T} \\
                                           & = \left( {1 \over m} \cdot (m, m, m, \dots, m, 0) \right) \cdot (m^n, m^{n-1}, m^{n-2}, \dots, m^0)^\mathsf{T} \\
                                           & = (1, \dots, 1, 0) \cdot (m^n, m^{n-1}, m^{n-2}, \dots, m^0)^\mathsf{T} \\
                                           & = \sum_{k=1}^{n} m^k = m \cdot {(m^n - 1) \over (m-1)}
\end{split}
\end{equation}

So, we have shown that $\mathbb{E}[\xi_{m,n}] = m \cdot {(m^n - 1) \over (m-1)}$.

Now, let's compute $\mathbb{E}[\xi_{m,n}^2]$:
\begin{equation}
\begin{split}
    \mathbb{E}[\xi_{m,n}^2] & = (0, \dots , 0, 1) \cdot W \cdot (\mathbb{I} + W) \cdot \left( (\mathbb{I} - W)^{-1} \right)^3 \cdot (1, 0, \dots , 0)^\mathsf{T} \\
                                           & = \left( 0, \dots, 0, {1 \over m}, 0 \right) \cdot (\mathbb{I} + W) \cdot \left( (\mathbb{I} - W)^{-1} \right)^2 \cdot (m^n, m^{n-1}, m^{n-2}, \dots, m^0)^\mathsf{T} \\
                                           & = \left( \left( 0, \dots, 0, {1 \over m}, 0 \right) \cdot (\mathbb{I} + W) \right) \cdot \left( (\mathbb{I} - W)^{-1} \right)^2 \cdot (m^n, m^{n-1}, m^{n-2}, \dots, m^0)^\mathsf{T} \\
                                           & = \left(0, \dots ,0, {1 \over m^2}, {1 \over m}, 0 \right) \cdot \left( (\mathbb{I} - W)^{-1} \right)^2 \cdot (m^n, m^{n-1}, m^{n-2}, \dots, m^0)^\mathsf{T} \\
                                           & = \left( \left(0, \dots ,0, {1 \over m^2}, {1 \over m}, 0 \right) \cdot (\mathbb{I} - W)^{-1} \right) \cdot (\mathbb{I} - W)^{-1} \cdot (m^n, m^{n-1}, m^{n-2}, \dots, m^0)^\mathsf{T} \\
                                           & = \left(2, \dots, 2, {2m - 1 \over m}, 0 \right) \cdot (\mathbb{I} - W)^{-1} \cdot (m^n, m^{n-1}, m^{n-2}, \dots, m^0)^\mathsf{T} \\
                                           & = (a_0, a_1, \dots, a_{n-1}, 0) \cdot (m^n, m^{n-1}, m^{n-2}, \dots, m^0)^\mathsf{T}
\end{split}
\end{equation}
Where $a_k = {2m \over m - 1} \cdot (m^n - m^k) - 1$ for every $k \in \{0, 1, \dots, n-1 \}$. \\

Thus:
\begin{equation}
\label{equation_e_xi_squared}
\begin{split}
    \mathbb{E}[\xi_{m,n}^2] & = \sum_{k=1}^{n} m^k \cdot a_{n-k} \\
                                           & = \sum_{k=1}^{n} m^k \cdot \left({2m \over m - 1} \cdot (m^n - m^{n - k}) - 1 \right) \\
                                           & = {m \over (m-1)^2} \cdot (2 \cdot m^{2n+1} -(2n+3) \cdot m^{n+1} + (2n+1) \cdot m^n + m - 1)
\end{split}
\end{equation}

From equations \eqref{equation_e_xi} and \eqref{equation_e_xi_squared} we can obtain the closed-form expression for variance:
\begin{equation}
\label{equation_variance_closed_form_derived}
\begin{split}
    \mathrm{Var}[\xi_{m,n}] & = \mathbb{E}[\xi_{m,n}^2] - \mathbb{E}[\xi_{m,n}]^2 \\
                                          & = {m \over (m-1)^2} \cdot (m^{2n+1} - (2n+1) \cdot m^{n+1} + (2n+1) \cdot m^n - 1)
\end{split}
\end{equation}

\section{Sum of the common sub-path lengths in the complete $m$-ary tree of depth $n$}
\label{paths_in_tree_section}

Given a complete $m$-ary rooted tree of depth $n$: $G_{m,n}=(V, E)$, we will find the closed-form expression for $S_{m,n} := \sum_{(a,b) \in V \times V} |\pi(a) \cap \pi(b)|$, where $\pi(v)$ is a set of edges that represents a path from the node $v \in V$ to the root node of the tree.

Consider some integer $d \in \{ 1, 2, \dots n \}$ which represents the length of a common sub-path. Let's define the subset $P_d \subset V \times V$, such that $\forall (a,b) \in P_d$: $|\pi(a) \cap \pi(b)| = d$. Then:
\begin{equation}
    S_{m,n} = \sum_{d=1}^n d \cdot |P_d|
\end{equation}

Consider a pair $(a, b) \in P_d$. As far as the length of the common sub-path between $\pi(a)$ and $\pi(b)$ is $d$, the deepest common node of these two paths (let's call it $c$) is located in the $d$-th level of the tree. There are $m^d$ nodes in the $d$-th level of the tree. Hence, the set $P_d$ can be represented as a union of $m^d$ disjoint subsets $Q_c$ (where $Q_c$ represents the set of pairs $(a,b)$ whose paths contain a deepest common node $c$, located in the $d$-th level of the tree):
\begin{equation}
\label{eq:number_of_sub_paths_with_length_d}
    |P_{d}| = m^d \cdot |Q_c|
\end{equation}

\begin{figure}
    \subfloat[Both nodes from the tuple $(a, b)$ are descendants of the node $c$.] 
    {{ \includegraphics[scale=0.35]{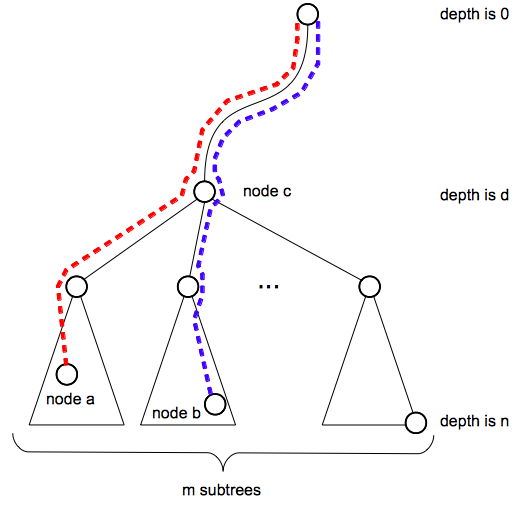} \label{fig:tree_case_1_figure} }}
    \qquad
    \subfloat[One node from the tuple $(a, b)$ is descendant of the node $c$, and the other node from the tuple is the node $c$.] 
    {{ \includegraphics[scale=0.4]{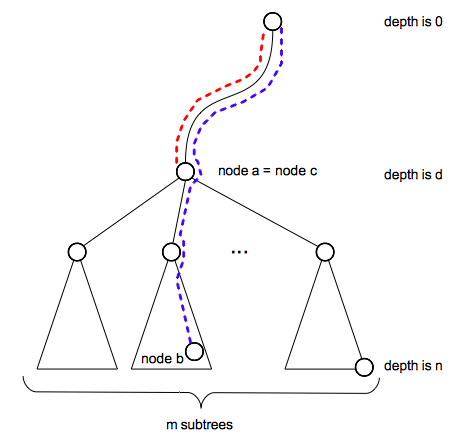} \label{fig:tree_case_2_figure} }}
    \caption{Two cases when the length of the common sub-path is $d$.}
\end{figure}

For computing $|Q_c|$ we need to consider the following three disjoint cases:
\begin{itemize}

    \item \textit{Case 1:} Nodes $a$ and $b$ are located in the different subtrees of $c$ (see the figure \ref{fig:tree_case_1_figure}). There are ${m \choose 2}$ possibilities to choose such subtrees. Each of these subtrees contains $\sum_{k=0}^{n-1-d} m^k = {m^{n-d} - 1 \over m - 1}$ nodes. The node $a$ is taken from one subtree, and the node $b$ is taken from the other subtree. The total number of tuples $(a, b)$ constructed in such case is: $2 \cdot {m \choose 2} \cdot \left( {m^{n-d} - 1 \over m - 1} \right)^2$.
    
    \item \textit{Case 2:} One of the nodes of the tuple is $c$, and the other node belongs to the subtree of $c$ (see the figure \ref{fig:tree_case_2_figure}). There are $m \cdot \sum_{k=0}^{n-1-d} m^k = {m^{n-d+1} - 1 \over m - 1} - 1$ possibilities to choose a descendant of $c$. The total number of tuples $(a, b)$ constructed in such case is: $2 \cdot \left( {m^{n-d+1} - 1 \over m - 1} - 1 \right)$.
    
    \item \textit{Case 3:} Both nodes in the tuple are $c$: $(a,b)=(c,c)$. There is only $1$ such case.
    
\end{itemize}

Thus, taking into account that the described cases are disjoint: 
\begin{equation}
\label{eq:number_of_sub_paths_with_length_d_and_common_vertex_c}
    |Q_c| = 2 \cdot {m \choose 2} \cdot \left( {m^{n-d} - 1 \over m - 1} \right)^2 +  2 \cdot \left( {m^{n-d+1} - 1 \over m - 1} - 1 \right) + 1
\end{equation}

From equations \eqref{eq:number_of_sub_paths_with_length_d} and \eqref{eq:number_of_sub_paths_with_length_d_and_common_vertex_c} we have:
\begin{equation}
    |P_d| = m^d \cdot \left( 2 \cdot {m \choose 2} \cdot \left( {m^{n-d} - 1 \over m - 1} \right)^2 +  2 \cdot \left( {m^{n-d+1} - 1 \over m - 1} - 1 \right) + 1 \right)
\end{equation}

Hence the sum of the common lengths over all 2-tuples of nodes of the complete $m$-ary tree of depth $n$ is:
\begin{equation}
\label{equation_equation_sum_common_subpaths_derived}
\begin{split}
    S_{m,n} & = \sum_{d=1}^{n} d \cdot m^d \cdot \left( 2 \cdot {m \choose 2} \cdot \left( {m^{n-d} - 1 \over m - 1} \right)^2 +  2 \cdot \left( {m^{n-d+1} - 1 \over m - 1} - 1 \right) + 1 \right) \\
                  & = \sum_{d=1}^{n} d \cdot m^d \cdot \left( { m^{2n-2d+1} - 1 \over m - 1} \right) \\
                  & = {m \over (m-1)^3} \cdot \left( m^{2n+1} - (2n+1) \cdot m^{n+1} + (2n+1) \cdot m^n - 1 \right)
\end{split}
\end{equation}

So, from the equations \eqref{equation_variance_closed_form_derived} and \eqref{equation_equation_sum_common_subpaths_derived} we see that $\mathrm{Var}[\xi_{m,n}] = (m - 1) \cdot S_{m,n}$.

In case when $m=2$ we have:
\begin{equation}
    \mathrm{Var}[\xi_{2,n}] = S_{2,n} = 4 \cdot 2^{2n} - (4n+2) \cdot 2^n - 2
\end{equation}

It is known that $S_{2,n}$ is described by the integer sequence A286778 \cite{oeis_A286778}.
But now, we have obtained a new interpretation for the integer sequence A286778: \textit{this sequence describes $\mathrm{Var}[\xi_{2,n}]$ -- a variance of the number of tosses of a fair coin until we see $n$ heads in a row}.

\section{Convergence of $\sum_{k=0}^{\infty} W^k$}
\label{infinite_series_convergence}

In this section we will show that the series $\sum_{k=0}^{\infty} W^k$ converges (taking into account the structure of the matrix $A$ from the equation \eqref{eq:matrix_A_structure} and that $W:={1 \over m} \cdot A$). Firstly, let's check how does the linear transformation $W$ act on some arbitrary column-vector  $\vec{w} = (f_k, f_{k-1}, f_{k-2}, \dots, f_{k-n})^\mathsf{T}$:
\begin{equation}
\begin{split}
     W \cdot \vec{w} = {1 \over m} \cdot A \cdot \vec{w} = {1 \over m} \cdot (f_{k+1}, f_k, f_{k-1}, f_{k-2}, \dots, f_{k-n + 1})^\mathsf{T}
\end{split}
\end{equation}
where $f_{k+1} = (m-1) \cdot \sum_{j=0}^{n-1} f_{k-j}$, and $m \geq 2$, $n \geq 1$. 

So, we see that the matrix $A$ encodes a linear recurrence relation: 
\begin{equation}
\label{eq:recurrence_relation_encoded_by_A}
    f_{k+1} = (m-1) \cdot \sum_{j=0}^{n-1} f_{k-j}
\end{equation}

In case if $m=2$ this is a recurrence relation for the Fibonacci $n$-Step Numbers \cite{generalized_fibonacci}. Below is a characteristic polynomial for the recurrence relation \eqref{eq:recurrence_relation_encoded_by_A}:
\begin{equation}
\label{eq:characteristic_polynomial}
\begin{split}
    x^n - (m-1) \cdot \sum_{j=0}^{n-1} x^j = 0
\end{split}
\end{equation}

As described in \cite{combinatorics_vilenkin}: when $k \to \infty$, the value of a $k$-th element of a linear recurrence relation is proportional to $|r|^k$, where $r$ is a root of the characteristic polynomial with the maximal absolute value among other roots. So, we need to analyse the roots of the characteristic polynomial \eqref{eq:characteristic_polynomial}.

In the section \ref{upper_bound_polynomial_section} we will show that the equation \eqref{eq:characteristic_polynomial} is equivalent to the equation $x^n \cdot (m - x) = m-1$, and we will show that for any $n$ the absolute values of roots of this polynomial are strictly less than $m$: $|r| < m$ (when $m=2$ this result is consistent to the result described in \cite{generalized_fibonacci}: the $k$-th Fibonacci $n$-Step Number is proportional to $r^k$, where $1 < r < 2$ is a solution of the equation $x^n \cdot (2-x) = 1$). \\

Consider a column-vector that has $1$ at the $i$-th position, and $0$ at all other positions: $\vec{v_i} = (0, 0, \dots, 0, 1, 0, \dots 0)^\mathsf{T} = (f_{n}, f_{n-1}, f_{n-2}, \dots, f_{0})^\mathsf{T}$. The expression $A^k \cdot \vec{v_i}$ allows to obtain a vector with the $(n+k)$-th, $(n+k-1)$-th (and so forth) elements of the recurrence relation \eqref{eq:recurrence_relation_encoded_by_A}: $(f_{n+k}, f_{n+k-1}, f_{n+k-2}, \dots, f_{k-n})^\mathsf{T}$, based on the given vector of initial values $\vec{v_i}$.

On the other hand, the expression $A^k \cdot \vec{v_i}$ equals to the $i$-th column of the matrix $A^k$. Hence, all elements of the matrix $A^k$ correspond to the items of the linear recurrence \eqref{eq:recurrence_relation_encoded_by_A}, but the values in each column are based on the different vectors of initial values $\vec{v_i}$. \\

Hence, all values in the matrix $A^k$ are non-negative and have an upper bound: $C \cdot |r|^{k+n}$ for some $|r| < m$, and where $C \in \mathbb{R}^{+}$ is some constant. Consequently, all values of the matrix $W^k$ are bound by ${1 \over m^k} \cdot C \cdot |r|^{k+n}$. Thus, the upper bound for the Frobenius norm of $W^k$ is:

\begin{equation}
    \Vert W^k \Vert_F \leq \sqrt{ \sum_{i=0}^n \sum_{j=0}^n \left( C \cdot {1 \over m^k} \cdot |r|^{n+k} \right)^2} = (n+1) \cdot C \cdot  {1 \over m^k} \cdot |r|^{n+k}
\end{equation}

Let's estimate the spectral radius of $W$ using the Gelfand's formula:
\begin{equation}
\begin{split}
    \rho(W) & = \lim_{k \to \infty}\Vert W^k \Vert_F^{1 / k} \\
                 & \leq \lim_{k \to \infty} \left( (n+1) \cdot C \cdot  {1 \over m^k} \cdot |r|^{n+k} \right)^{1 / k} \\
                 & = \lim_{k \to \infty} {|r| \over m} \cdot \left( (n+1) \cdot C \cdot  |r|^n \right)^{1/k} \\
                 & = {|r| \over m} < 1
\end{split}
\end{equation}

As far as $\rho(W) < 1$ the series $\sum_{k=0}^{\infty} W^k$ converges to $(\mathbb{I} - W)^{-1}$ and is known as Neumann series.

\section{Upper bound on the absolute values of roots of the characteristic polynomial}
\label{upper_bound_polynomial_section}


Consider the polynomial \eqref{eq:characteristic_polynomial}: $x^n - (m-1) \cdot \sum_{j=0}^{n-1} x^j = 0$ where $m \geq 2$ and $n \geq 1$. The Cauchy's bound gives a non-strict bound for the absolute values of roots of this polynomial: the values are less or equal to $m$. Let's show that this bound is strict (that the absolute values of roots are strictly less than $m$).

As far as $1$ is not a root of this polynomial, we can rewrite this equation as follows:
\begin{equation}
\begin{split}
    x^n - (m-1) \cdot {x^n - 1 \over x - 1}   & = 0 \\
    \iff x^n \cdot (m - x) &= m-1
\end{split}
\end{equation}

We need to consider the roots of the polynomial: $x^n \cdot (m - x) - (m-1) = 0$.
For the sake of contradiction let's assume, that $x =m \cdot e^{i \phi}$ (for some $\phi \in \mathbb{R}$) is a root of this polynomial.
After doing the substitution we obtain:
\begin{equation}
\begin{split}
          m^n \cdot e^{i n \phi} \cdot (m - m \cdot e^{i \phi}) &= m - 1 \\
    \iff  e^{i n \phi} - e^{i (n+1) \phi} &= {m - 1 \over m^{n+1}}
\end{split}
\end{equation}

Using the trigonometric form of complex numbers we can rewrite the latter equality as follows:
\begin{equation}
    \left( cos(n\phi) - cos((n+1)\phi) \right) - i \cdot \left( sin(n\phi) - sin((n + 1)\phi) \right) = {m - 1 \over m^{n+1}}
\end{equation}

As far as the imaginary part on the left hand side is $0$ we have a following system:
\begin{equation}
    \left\{
    \begin{aligned}
        sin(n\phi) - sin((n+1)\phi) &= 0 \\ 
        cos(n\phi) - cos((n+1)\phi) &= {m - 1 \over m^{n+1}}
    \end{aligned}\right.
\end{equation}

Using the sum-to-product trigonometric identities we can rewrite the system as follows:
\begin{equation}
    \left\{
    \begin{aligned}
         2 \cdot cos \left({2n + 1 \over 2} \phi \right) \cdot sin \left( - {\phi \over 2} \right) &= 0 \\ 
        -2 \cdot sin \left({2n + 1 \over 2} \phi \right) \cdot sin \left( - {\phi \over 2} \right) &= {m - 1 \over m^{n+1}}
    \end{aligned}\right.
\end{equation}

Looking at the first equation: $2 \cdot cos \left({2n + 1 \over 2} \phi \right) \cdot sin \left( - {\phi \over 2} \right)$ we conclude that some of its multipliers should be $0$. So, we have two cases:
\begin{itemize}

    \item \textit{Case 1}: $sin \left( - {\phi \over 2} \right) = 0$. In this case the second equation shows a contradiction: $-2 \cdot sin \left({2n + 1 \over 2} \phi \right) \cdot sin \left( - {\phi \over 2} \right) = 0 \neq {m - 1 \over m^{n+1}}$

    \item \textit{Case 2}: $cos \left({2n + 1 \over 2} \phi \right) = 0$. In this case we can see that: $\phi = {2k + 1 \over 2n + 1} \pi$ for $k \in \mathbb{Z}$. Let's substitute $\phi$ into the second equation:
        \begin{equation}
        \begin{split}
                -2 \cdot sin \left({2n + 1 \over 2} \cdot {2k + 1 \over 2n + 1} \pi \right) \cdot sin \left( - {1 \over 2} \cdot {2k + 1 \over 2n + 1} \pi \right) &= {m - 1 \over m^{n+1}} \\
                \iff sin \left({2k + 1 \over 2 \cdot (2n + 1)} \pi \right) &= {m - 1 \over 2 \cdot m^{n+1}}
        \end{split}
        \end{equation}
    So, from the second equation we have obtained the equation $sin(a \cdot \pi) = b$, where $a$ and $b$ are rational numbers ($a={2k + 1 \over 2 \cdot (2n + 1)}$ and $b={m - 1 \over m^{n+1}}$). Furthermore, in case if $m \geq 2, n \geq 1$ we see that $b \not \in \{ 0, \pm 1, \pm {1 \over 2}\}$. In this case, we have a contradiction with the Niven's Theorem \cite{nivens_theorem} (that states, that if $sin(a \cdot \pi) = b$ and $a, b \in \mathbb{Q}$, then the sine takes only the values $0, \pm 1, \pm {1 \over 2}$).
    
\end{itemize}

So, we see that the original assumption (that $m \cdot e^{i \phi}$ is a root of the polynomial) leads to the contradiction. Consequently, $m \cdot e^{i \phi}$ can't be a root of the polynomial.
In combination with the Cauchy's bound we have a strict bound on the absolute values of roots: the roots are strictly less than $m$.

\section{Matrix products properties}
\label{matrix_products_properties_section}

The product of matrices $W \cdot (\mathbb{I} - W)^{-1}$ is commutative: 
\begin{equation}
\label{equation_w_times_one_minus_w_inv_commutative}
\begin{split}
    W \cdot (\mathbb{I} - W)^{-1}                     & = (\mathbb{I} - (\mathbb{I} - W)) \cdot (\mathbb{I} - W)^{-1} \\
                                             	         	    & = (\mathbb{I} - W)^{-1} - (\mathbb{I} - W) \cdot (\mathbb{I} - W)^{-1} \\
                                         			            & = (\mathbb{I} - W)^{-1} - (\mathbb{I} - W)^{-1} \cdot (\mathbb{I} - W) \\
		                                                      & = (\mathbb{I} - W)^{-1} \cdot (\mathbb{I} - (\mathbb{I} - W)) \\
                    		                                     & = (\mathbb{I} - W)^{-1} \cdot W
\end{split}
\end{equation}

The product of matrices $W \cdot (\mathbb{I} + W)$ is commutative: 
\begin{equation}
\begin{split}
    W \cdot (\mathbb{I} + W) = W + W^2 = (\mathbb{I} + W) \cdot W
\end{split}
\end{equation}

Taking into account the commutativity property \eqref{equation_w_times_one_minus_w_inv_commutative} we can show that the products of matrices $(\mathbb{I} + W) \cdot (\mathbb{I} - W)^{-1}$ is also commutative:
\begin{equation}
\begin{split}
    (\mathbb{I} + W        ) \cdot (\mathbb{I} - W)^{-1}              & = (\mathbb{I} - W)^{-1} + W \cdot (\mathbb{I} - W)^{-1} \\
                                                                                               & = (\mathbb{I} - W)^{-1} + (\mathbb{I} - W)^{-1} \cdot W \\
                                                                                               & = (\mathbb{I} - W)^{-1} \cdot (\mathbb{I} + W)
\end{split}
\end{equation}

\end{document}